\documentclass[final]{article}
\usepackage {latexsym} \usepackage {amsmath} \usepackage {amsfonts} 
\usepackage {epsf}
\usepackage {a4}
\usepackage {amssymb}
\usepackage {listings}
\usepackage {enumerate}
\usepackage {graphicx}
\usepackage {QED}
\usepackage {url}
\def\code#1{\lstinline!#1!}
\def\zet{{\mathbb Z}}

\def\MC{\mathop{\rm MAXCUT}\nolimits}

\def\Nese{Ne\v set\v ril}
\def\sym{\mathbin{\Delta}}
\def\fsets{{{\cal P}([4])}}

\def\Pm{{\rm Parent\_menu}}
\def\lesseq{\preccurlyeq}

\def\hom{\mathrel{\xrightarrow{hom}}}

\def\cost{\mathrm{cost}}

\newtheorem{theorem}{Theorem}[section]
\newtheorem{lemma}[theorem]{Lemma}
\newtheorem{remark}[theorem]{Remark}

\newtheorem{proposition}[theorem]{Proposition}

\newtheorem{conjecture}[theorem]{Conjecture}
\newtheorem{observation}[theorem]{Observation}
\newtheorem{problem}[theorem]{Problem}

\newtheorem{claim}[theorem]{Claim}

\newenvironment{proofof}{\par\medskip\proofof}{\qed}
\def\proofof #1{\noindent{\bf Proof of #1:}}
\def\?#1{\par{\bf ??? #1 !!!}\par}
\begin{document}
\title{High-girth cubic graphs are homomorphic to the Clebsch graph}
\author{
  Matt DeVos
\thanks{Department of Mathematics, Simon Fraser University, Burnaby, B.C. V5A 1S6, Canada.
   Email: {\tt mdevos@sfu.ca}}
  \and
	Robert \v S\'amal
\thanks{Department of Applied Mathematics and Institute for Theoretical Computer Science (ITI),
 Charles University, 
 Malostransk\'e n\'am\v est\'\i{}~25, 118~00, Prague, Czech Republic.
 Email: {\tt samal@kam.mff.cuni.cz} }
\thanks{Institute for Theoretical Computer Science is supported as project
 1M0545 by Ministry of Education of the Czech Republic.}
}
\date{}
\maketitle
\leftline{{\bf Keywords:} max-cut, cut-continuous mappings, circular chromatic number}
\leftline{{\bf MSC:} 05C15}
\begin{abstract}
We give a (computer assisted) proof that the edges of every graph with maximum
degree 3 and girth at least 17 may be 5-colored (possibly improperly) so that
the complement of each color class is bipartite.  Equivalently, every such
graph admits a homomorphism to
the Clebsch graph (Fig.~\ref{fig:PetCl}).  

Hopkins and Staton \cite{HS} and Bondy and Locke \cite{BL}
proved that every (sub)cu\-bic graph of girth at least 4 has an edge-cut
containing at least $\frac{4}{5}$ of the edges.  The existence of such
an edge-cut follows immediately from the existence of a
5-edge-coloring as described above, so our theorem may be viewed as a
coloring extension of their result (under a stronger girth
assumption).  

Every graph which has a homomorphism to a cycle of length five has 
an above-described 5-edge-coloring; hence our theorem may also be viewed
as a weak version of \Nese's Pentagon Problem (which asks whether every 
cubic graph of sufficiently high girth is homomorphic to~$C_5$). 
\end{abstract}

\section{Introduction}\label{sec:intro}

Throughout the paper all graphs are assumed to be finite, undirected
and simple.  For any positive integer $n$, we let $C_n$ denote the
cycle of length $n$, and $K_n$ denote the complete graph on $n$
vertices.  If $G$ is a graph and $U \subseteq V(G)$, we put $\delta(U)
= \{ uv \in E(G) : \mbox{$u \in U$ and $v \not\in U$} \}$, and we call
any subset of edges of this form a \emph{cut}.  The maximum size of a cut
of $G$, denoted $\MC (G) = \max_{U \subseteq V} |\delta(U)|$ is a
parameter which has received great attention.  Next, we normalize and
define
$$
    b(G) = \frac {\MC (G)}{ |E(G)| } \,.
$$
Determining $b(G)$ (or equivalently $\MC(G)$) for a given graph~$G$ is
known to be NP-complete, so it is natural to seek lower bounds.  It is
an easy exercise to show that $b(G) \ge 1/2$ for any graph~$G$ and
$b(G) \ge 2/3$ whenever $G$ is cubic (that is 3-regular).  The former
inequality is almost attained by a large complete graph, the latter is
attained for $G = K_4$: any triangle contains at most two edges from
any bipartite subgraph, and each edge of~$K_4$ is in two triangles.  
This suggests that triangles play a special role,
and raises the question of improving this bound for cubic graphs of
higher girth. In the 1980's, several authors independently considered
this problem \cite{BL,HS,Zy}, the strongest results being
\begin{itemize}
  \item $b(G) \ge 4/5$ for $G$ with maximum degree 3 and
     no triangle \cite{BL}
  \item $b(G) \ge 6/7 - o(1)$ for cubic $G$ with girth tending to
     infinity \cite{Zy} \end{itemize}
On the other hand, cubic graphs exist with arbitrarily high girth and
satisfying $b(G) < 0.9351$. 
This result was announced by McKay in a conference~\cite{McK-proc} but did not 
appear in print. It is, however, rather straightforward to prove it (with a 
worse constant) by considering random cubic graphs, see~\cite{Wormald-survey} 
for a nice survey. In thesis of Jan Hladk\'y~\cite{Hladky} the constant $0.9351$ 
is recovered.

Define a set of edges $C$ from a graph $G$ to be a \emph{cut
complement} if $C = E(G) \setminus \delta(U)$ for some $U \subseteq
V(G)$.  Then the problem of finding a cut of maximum size is exactly
equivalent to that of finding a cut complement of minimum size.  A
natural relative of this is the problem of finding many disjoint cut
complements.  Indeed, packing cut complements may be viewed as a 
coloring version of the maximum cut problem.  

There are a variety of interesting properties which are equivalent to
the existence of $2k+1$ disjoint cut complements, so after a handful
of definitions we will state a proposition which reveals some of these
equivalences.  This proposition is well known, but we have provided a
proof of it in Section~\ref{sec:equiv} for the sake of completeness.  For
every positive integer $n$, we let $Q_n$ denote the \emph{$n$-dimensional
cube}, so the vertex set of $Q_n$ is the set of all binary vectors of
length $n$, and two such vertices are adjacent if they differ in a
single coordinate. 
The \emph{$n$-dimensional projective cube},\footnote{sometimes called folded cube}
denoted $PQ_n$, is the simple graph obtained from the $(n+1)$-dimensional
cube~$Q_{n+1}$ by identifying pairs of antipodal vertices (vertices that
differ in all coordinates). Equivalently, the projective cube $PQ_n$ can be 
described as a Cayley graph, see Section~\ref{sec:equiv}.
If $G$, $H$ are graphs, a
\emph{homomorphism} from $G$ to $H$ is a mapping $f: V(G) \rightarrow
V(H)$ with the property that $f(u)f(v)$ is an edge of $H$ whenever
$uv$ is an edge of $G$.  
When there exists a homomorphism from $G$ to $H$, we say that
\emph{$G$ is homomorphic to $H$} and write $G \to H$. 
We need yet another concept, introduced in \cite{DNR}:
A mapping $g : E(G) \rightarrow E(H)$ is
\emph{cut-continuous} if the preimage of every cut is a cut.  Now we
are ready to state the relevant equivalences.

\begin{proposition} \label{equiv_prop}
For every graph $G$ and nonnegative integer $k$, the following
properties are equivalent.  
\begin{enumerate}[(1)]
  \item There exist $2k$ pairwise disjoint cut complements.
  \item There exist $2k+1$ pairwise disjoint cut complements with union $E(G)$.
  \item $G$ has a homomorphism to $PQ_{2k}$.
  \item $G$ has a cut-continuous mapping to $C_{2k+1}$.  
\end{enumerate}
\end{proposition}

\bigskip

Perhaps the most interesting conjecture concerning the packing of cut
comple\-ments---or equivalently homomorphisms to projective cubes---is
the following conjectured generalization of the Four Color Theorem.
Although not immediately obvious, this is equivalent to Seymour's \cite{seymour-r} 
conjecture on $r$-edge-coloring of planar $r$-graphs (when $r$ is odd).

\begin{conjecture}[Seymour]
Every planar graph with all odd cycles of length at least~$2k+1$ has a
homomorphism to $PQ_{2k}$.  
\end{conjecture}

\begin{figure}
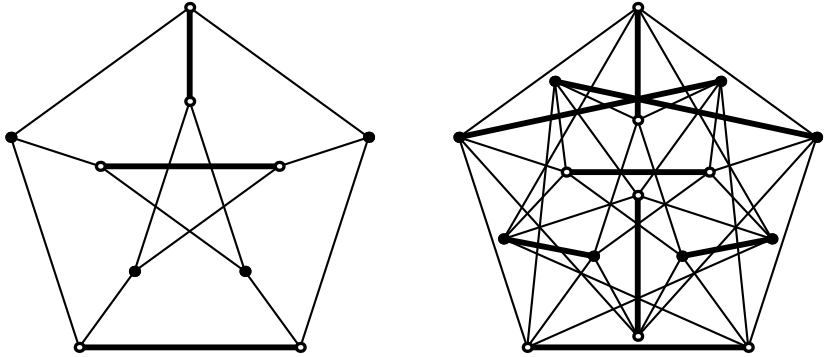

\centerline{\includegraphics{Graphs.0} \hfil \includegraphics{Graphs.1}}
\caption {Petersen and Clebsch graph with one cut complement
emphasized, the respective bipartition of the vertex set is
depicted, too. The other four cut complements are obtained by 
rotation.}
\label{fig:PetCl}
\end{figure}

Since the graph $PQ_2$ is isomorphic to $K_4$, the $k=1$ case of this
conjecture is equivalent to the Four Color Theorem.  The $k=2$ case of
this conjecture concerns homomorphisms to the graph $PQ_4$ which is
also known as the Clebsch graph (see Figure~\ref{fig:PetCl}).  This
case was resolved in the affirmative by Naserasr~\cite{Naserasr-thesis}
who deduced it from a theorem of Guenin~\cite{Guenin-joins}.  

The following theorem is the main result of this paper; it shows that graphs of
maximum degree three without short cycles also have homomorphisms to
$PQ_4$.  The \emph{girth} of a graph is the length of its shortest 
cycle, or $\infty$ if none exists.

\begin{theorem} \label{main}
Every graph of maximum degree~3 and girth at least~$17$ is homomorphic
to $PQ_4$ (also known as the Clebsch graph), or equivalently has 5
disjoint cut complements.  Furthermore, there is a linear time
algorithm which computes the homomorphism and the cut complements.
\end{theorem}

Clearly no graph with a triangle can map homomorphically
to the triangle-free Clebsch graph (equivalently, have $5$ disjoint cut
complements), but we believe this to be the only obstruction for cubic
graphs. We highlight this and one other question we have been unable
to resolve below.  

\begin{conjecture}[\cite{rs-thesis}]
Every triangle-free cubic graph is homomorphic to $PQ_4$. 
\end{conjecture}

\begin{problem}
What is the largest integer $k$ with the property that all cubic
graphs of sufficiently high girth are homomorphic to $PQ_{2k}$?
\end{problem}

As we mentioned before, there are high-girth cubic graphs $G$
with $b(G) < 0.94$. Such graphs do not admit homomorphism to
$PQ_{2k}$ for any $k \ge 8$ (Proposition~\ref{equiv_prop} (2)), 
so there is indeed some largest integer~$k$ in the above problem. 
At present, we know only that $2 \le k \le 7$.  

\bigskip

Another topic of interest for cubic graphs of high girth is circular
chromatic number, a parameter we now pause to define.  For any graph
$G$, we let $G^{\ge k}$ denote the simple graph with vertex set $V(G)$
and two nodes adjacent if they have distance at least~$k$ in~$G$.  The
\emph{circular chromatic number} of $G$, is $\chi_c(G) = \inf \{
\frac{n}{k} : \mbox{$G$ has a homomorphism to $C_n^{\ge k}$} \}$.
Every graph satisfies $\lceil \chi_c(G) \rceil = \chi(G)$ so the
circular chromatic number is a refinement of the usual notion of
chromatic number.  The following curious conjecture asserts that cubic
graphs of sufficiently high girth have circular chromatic number $\le
\frac{5}{2}$ (since $C_{2k+1}$ and  $C_{2k+1}^{\ge k}$ are isomorphic).

\begin{conjecture}[\Nese's Pentagon Conjecture \cite{Nes-Aspects}]   \label{hompentagon}
If $G$ is a cubic graph of sufficiently high girth
then $G$ is homomorphic to~$C_5$.
\end{conjecture}

It is an easy consequence of Brook's Theorem that the above conjecture
holds with~$C_3$ in place of~$C_5$ (i.e., every cubic graph of girth at
least~$4$ is 3-colorable). On the other hand, it is known that the
conjecture is false if we replace~$C_5$ by $C_7$ \cite{Hatami-C7},
consequently it is false if we replace $C_5$ by any $C_n$ for
odd $n \ge 7$. (Earlier, this result was proved for 
$n\ge 11$ \cite{KNS} and $n \ge 9$ \cite{WW}.)

An important extension of Conjecture~\ref{hompentagon} is the problem
to determine the infimum of real numbers $r$ with the property that every cubic
graph of sufficiently high girth has circular chromatic number $\le r$.  The
above results show that this infimum must lie in the interval
$[\frac{7}{3},3]$, but this is the extent 
of our knowledge.  
It is tempting to try to use the fact that girth $\ge 17$ cubic graphs
map to the Clebsch graph and girth $\ge 4$ cubic graphs map to $C_3$
to improve the upper bound, but the circular chromatic numbers of
$C_3$, the Clebsch graph, and their direct product are all at
least three,\footnote{The only nontrivial case is the product
$PQ_4\times{}K_3$. By a theorem of~\cite{GL-unique} this graph is
uniquely 3-colorable; consequently $\chi_c(PQ_4\times{}K_3) = 3$.} 
so no such improvement can be made. 
Neither were we able to use our result to improve upper bounds on
fractional chromatic numbers of cubic graphs. This is conjectured to
be at most~$14/5$ for triangle-free cubic graphs 
(Heckmann and Thomas~\cite{HT-independenceratio}), and proved to be
at most~$3-3/64$ (Hatami and Zhu~\cite{HatamiZhu}).

It is easy to prove directly that Conjecture~\ref{hompentagon}, if
true, implies Theorem~\ref{main} (perhaps with a stronger assumption
on the girth).  This follows from part (4) of Proposition
\ref{equiv_prop} and the following easy observation.

\begin{observation}  \label{homocc}
If there is a homomorphism from $G$ to $H$, then there is a cut-continuous mapping from 
$G$ to $H$.  
\end{observation}

\begin{proof}
Let $f: V(G) \to V(H)$ be a homomorphism and define the mapping $f^\sharp: E(G) \to E(H)$ 
by the rule $f^\sharp(uv) = f(u)f(v)$.  If $S=\delta(U)$ is a cut in~$H$,  
then $(f^\sharp)^{-1}(S) = \delta(f^{-1}(U))$, which is also a cut.
\end{proof}

The relationship between homomorphisms and cut-continuous maps is
studied in greater detail in~\cite{NS-TT1} and~\cite{NS-TT2} where it
is shown that, perhaps surprisingly, existence of a cut-continuous
mapping from~$G$ to~$H$ frequently implies the existence of a
homomorphism from~$G$ to~$H$.  Unfortunately, it does not appear
likely that these techniques can be used to extend the main theorem of
this paper to attain Conjecture~\ref{hompentagon}.

We finish the introduction with another conjecture due to Ne\v set\v ril 
(personal communication) concerning the existence of homomorphisms for cubic 
graphs of high girth.

\begin{conjecture}   \label{sparseconj}
For every integer $k$ there is a graph~$H$ of girth at least~$k$ and
an integer~$N$, such that for every cubic graph~$G$ of girth at
least~$N$ we have
$$
        G \hom H \,.
$$
\end{conjecture}

Let us note that if we replace ``girth'' by ``odd-girth'', than the
result is true, by a result of \cite{HH}, or in a greater 
generality \cite{dMN}.
In~\cite{HH} they also give a simple and explicit construction 
of such graph $H$ when $k=4$; a corollary of our main theorem 
is that one can let $H$ be the Clebsch graph (quite a bit
smaller than the graph constructed in~\cite{HH}) in this case.

\section{The Proof}\label{sec:pf}

The goal of this section is to prove the main theorem.  We begin with a lemma which
reduces our task to cubic graphs.

\begin{lemma}   \label{reduction}
If Theorem~\ref{main} holds for every cubic~$G$ then it holds
for every subcubic~$G$, too.
\end{lemma}

\begin{proof}
Let $G$ be a subcubic graph of girth at least~17. We will find a cubic
graph~$G'$ such that girth of~$G'$ is at least~17 and $G' \supseteq G$.
The lemma then follows, as restriction of any homomorphism
$G' \hom PQ_4$ to~$V(G)$ is the desired homomorphism $G \hom PQ_4$.

To construct~$G'$, put $r = \sum_{v \in V(G)} (3 - \deg(v))$. 
Let $H$ be an $r$-regular graph of girth at least~17 (it is well known that
such graphs exist, see, e.g., \cite{Biggs-cubic} for a nice survey). We
take $|V(H)|$ copies of~$G$.  For every edge~$uv$ of~$H$ we choose two
vertices of degree less than~3, one from a copy of~$G$ corresponding
to each of~$u$ and~$v$; then we connect these by an edge. Clearly,
this process will lead to a cubic graph containing~$G$ and with girth
at least the minimum of girths of~$G$ and~$H$.
\end{proof}

\bigskip

\paragraph{Proof outline: }
To show that cubic graphs of girth $\ge 17$ have homomorphism to the
Clebsch graph, we shall use property (1) from Proposition~\ref{equiv_prop} --- 
that is, we try to find a 4-tuple of pairwise disjoint cut complements.
A natural way to do so is to consider any 4-tuple of cut complements
and then make them as disjoint as possible. To say this precisely
we introduce several terms to describe the tuples of cut complements
and to measure ``how disjoint'' they are.

A \emph{labeling} of a graph $G$ is a four-tuple
$X = (X_1,X_2,X_3,X_4)$ so that each $X_i$ is a subset of $E(G)$.  
We call a labeling $X$ a \emph{cut labeling} if every $X_i$ is a cut,
and a \emph{cut complement labeling} if every $X_i$ is a cut complement.  If 
$X_i \cap X_j = \emptyset$ whenever $1 \le i < j \le 4$ we say that the 
labeling is \emph{wonderful}.  

Define function $a: \{0,1,\ldots,4\} \rightarrow \zet$ by 
$a(0) = 0$, $a(1) = 1$, $a(2) = 10$, $a(3) = 40$, and $a(4) = 1000$.
Now, for any labeling $X$, we define the \emph{label} of an edge $e$
(with respect to $X$) to be $l_X(e) = \{ i \in \{1,2,3,4\} : e \in X_i \}$, 
the \emph{weight} of $e$ to be $w_X(e) = |l_X(e)|$, and the
\emph{cost} of $e$ to be $\cost_X(e) = a(w_X(e))$.  Finally, we define
the \emph{cost} of $X$ to be $\cost(X) = \sum_{e \in E(G)} \cost_X(e)$.

The structure of our proof is quite simple: we prove that any cut
complement labeling of minimum cost in a cubic graph of girth $\ge 17$
is wonderful.  To show that such a labeling is wonderful, we shall
assume it is not, and then make a small local change to improve the
cost---thus obtaining a contradiction. (This also leads to a 
linear time algorithm. Confirming the outline above, each step of the
algorithm is making the four cut complements more disjoint, in the 
sense that it decreases the cost defined in the previous paragraph.)

The observation below will be used to make our local changes.  
For any sets $A,B$ we let $A \sym B = (A \setminus B) \cup (B \setminus A)$  
be the symmetric difference.
If $X = (X_1,\ldots,X_4)$ and $Y = (Y_1,\ldots,Y_4)$ are
labelings, then we let $X \sym Y = (X_1 \sym Y_1, \ldots, X_4 \sym Y_4)$.  

\begin{observation}   \label{1cut}
If $C$ is a cut and $D$ is a cut complement, then $C \sym D$ is a cut
complement.  Similarly, if $X$ is a cut complement labeling and $Y$ is
a cut labeling, then $X \sym Y$ is a cut complement labeling.
\end{observation}

\begin{proof}
Let $C = \delta(U)$ and $D = E(G) \setminus \delta(V)$.  Then 
$C \sym D = E(G) \setminus \bigl(\delta(U) \sym \delta(V)\bigr) 
          = E(G) \setminus \delta(U \sym V)$
so it is a cut complement. For labelings we consider each coordinate
separately.
\end{proof}

The graphs we consider have high girth, so they `locally are trees'.
Our proof will exploit this by using the above
observation to make changes on a tree.  

For example, consider the tree on Figure~\ref{fig:badlab} (on the top). 
This is supposed to be a part of a large cubic graph~$G$
with a corresponding part of a cut complement
labeling of~$G$. The dashed lines indicate two cuts of~$G$:
$Y_2$ (indicated by $\{2\}$) and $Y_3$ (indicated by $\{3\}$). 
Putting $Y_1 = Y_4 = \emptyset$, we get a cut labeling $(Y_1, Y_2, Y_3, Y_4)$
of~$G$. By Observation~\ref{1cut}, $X \sym Y$ is a cut complement labeling. 
It is easy to verify that $X \sym Y$ has lower cost that $X$: the number of 
weight~1 edges (edges contained in exactly one of the cut complements)
decreases by~1, the number of edges of weight 2, 3, and 4 is not changed.

Perhaps surprisingly, it is possible to reach a wonderful labeling
(a 4-tuple of disjoint cut complements) by a series of such \emph{local} operations.
We need, however, to get a bit more precise to describe how operations with trees
correspond to local operations with graphs. 

\begin{figure}[ht]
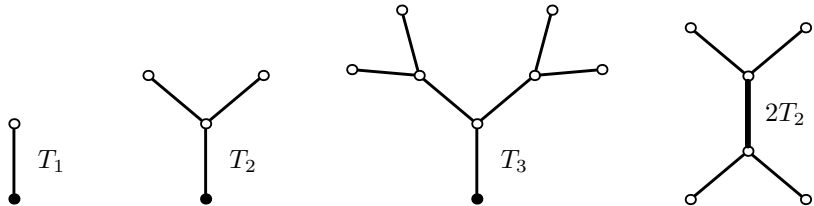

\centerline{
  \hfill
  \includegraphics{tree.5}
  \hfill
  \includegraphics{tree.6}
  \hfill
  \includegraphics{tree.7}
  \hfill
  \includegraphics{tree.8}
  \hfill
}
\caption {Illustration of definitions, root vertex/central edge are emphasized.}
\label{fig:defT_i}
\end{figure}

To this end, we
introduce a family of rooted trees (see Figure~\ref{fig:defT_i}).
Let $T_i$ denote a rooted
tree of ``depth $i$'' in which all vertices have degrees 1 and 3, and
the root vertex, denoted $r$, has degree 1.  Explicitly, we let $T_1$
be an edge (with one end being the root).  Having defined $T_i$, we
form $T_{i+1}$ by joining two copies of~$T_i$ by identifying their
root vertices and then connecting this common vertex to a new vertex,
which will be the new root.  The unique edge incident with the root we
shall call the \emph{root edge}.  We let $2T_i$ denote the tree
obtained from two copies of~$T_i$ by identifying their root edges in
the opposite direction (the resulting edge will be called the
\emph{central} edge of~$2T_i$).  A vertex of $T_i$ or $2T_i$ is
\emph{interior} if either it has degree 3, or it is the root of $T_i$.

A cut $C$ of $T_i$ or $2T_i$ is called \emph{internal} if
$C = \delta(Z)$ for some set $Z$ of interior vertices.  
A cut labeling $X$ is \emph{internal} if $X$ is a 4-tuple
of internal cuts.
(As illustrated above on the example from Figure~\ref{fig:badlab}, internal
cuts of a tree~$T$ correspond to ``normal'' cuts in a graph that
contains $T$ as a subgraph, possibly with some leafs identified. This is
utilized later, in the proof of Theorem~\ref{main}.)

Now we are ready to state and prove a lemma that forms the first step of the proof:
it will be used to show that any cut complement labeling of minimum cost has 
no edges of weight $> 2$.  

\begin{lemma}   \label{noweight34}
Let $X$ be a labeling of the tree $2T_2$ and assume that the weight of
the central edge is $>2$.  Then there exists an internal cut labeling
$Y$ of $2T_2$ so that $\cost(X \sym Y) < \cost(X)$.
\end{lemma}

(Note that we will actually prove this for $T_2$ in place of $2T_2$.
This version, however, corresponds better to Lemma~\ref{noweight2}.)

\begin{proof}
Let $e$ be the central edge, let $x$ be a vertex incident with $e$, let
$f,g$ be the other edges incident with $x$, and let $A = l_X(e)$, $B = l_X(f)$,
and $C = l_X(g)$.  We will construct an internal cut labeling 
$Y = (\delta(Z_1), \ldots, \delta(Z_4))$ 
(where each $Z_i$ is either $\emptyset$ or $\{x\}$)
so that $\cost(X \sym Y) < \cost(X)$.  
For convenience, we shall say that we \emph{switch} a set 
$I \subseteq \{1,2,3,4\}$ if we set $Z_i = \{x\}$ if $i \in I$ and
$Z_i = \emptyset$ otherwise.  

If $S=A \cap B \cap C$ is nonempty then we may switch $S$, thereby
reducing the cost of each of $e$, $f$, $g$.  Hence we may suppose $S$~is empty.

{\bf Case 1. $|A| = 4$:}
If $B = C = \emptyset$ then we switch $\{1\}$ decreasing the cost
from~$a(4)$ to~$a(3) + 2a(1)$. Otherwise we switch $B\cup C$; this
leads to a label $\{1,2,3,4\} \setminus (B \cup C)$ on~$e$, $C$ on~$f$ and $B$
on~$g$, reducing the cost again.

{\bf Case 2. $|A| = 3$:} We may suppose $A = \{1,2,3\}$ and 
$|A| \ge |B| \ge |C|$.
Moreover, $|C| < 3$ for otherwise $A \cap B \cap C$ is nonempty.
If $A$ and~$B$ have a common element, then we switch it.  
This changes the weights of edges in~$T$ from $3$, $|B|$, $|C|$ to $2$, $|B|-1$, 
$|C|+1$ and as $|C| < 3$, this
is an improvement in the total cost. It remains to consider the 
cases when both~$B$ and~$C$ are subsets of $\{4\}$.  
In each of these cases we switch $\{1\}$, this reduces the cost from 
at least~$a(3)$ to at most~$3a(2)$.  
\end{proof}

The next lemma, which provides the second step of the proof, is analogous 
to the previous one, but is considerably more complicated to prove. 

\begin{lemma}   \label{noweight2}
Let $X$ be a labeling of the tree $2T_9$ and assume that every edge has weight $\le 2$ 
and that the central edge has weight exactly 2.  Then there exists an internal cut labeling 
$Y$ of $2T_9$ so that $\cost(X \sym Y) < \cost(X)$.
\end{lemma}

Before discussing the proof of this lemma we shall use it to prove the main theorem.

\begin{proof}{\bf (of Theorem~\ref{main})}
It follows from Lemma~\ref{reduction} and Proposition~\ref{equiv_prop}
that it suffices to prove that all 
cubic graphs with girth at least~$17$ have wonderful cut complement
labelings. Let $G$ be such graph and let $X$ be a cut complement
labeling of $G$ of minimum cost.  It follows immediately from
Lemma~\ref{noweight34} that every edge of $G$ has weight $\le 2$.
Suppose there is an edge $e$ of weight $2$.  Then it follows from our
assumption on the girth that $G$ contains a subgraph isomorphic to
$2T_9$ (possibly with some of the leaf vertices identified) where $e$
is the central edge. Now Lemma~\ref{noweight2} gives us 
an \emph{internal} cut labeling~$Y$ of~$2T_9$ (hence a cut labeling of~$G$)
such that $\cost (X \sym Y) < \cost (X)$.
This contradiction shows that $X$ is wonderful, and completes the proof 
of the first part.

Next we give a short description of a linear-time algorithm that finds the 
partition. We start with a cut complement labeling $X=(E(G),E(G),E(G),E(G))$.  
Then we repeatedly pick a bad edge~$e$---that is an edge for which~$w(e) > 1$.
By Lemma~\ref{noweight34} and~\ref{noweight2} we can decrease
the total cost by moving from $X$ to $X \sym Y$ where $Y$ is a cut labeling 
that contains only edges at distance at most 8 from~$e$. We can therefore find the cut
labeling in constant time (e.g., by brute force if we do not try to minimize
the constant)---we only have to use an efficient representation of the
graph, namely a list of edges, list of vertices, and pointers between
the adjacent objects. As the cost of the starting coloring is 
$a(4)\cdot |E(G)|$ and at each step the decrease is at least by~1,
it remains to handle the operation ``pick a bad edge'' in constant time.
For this, we maintain a linked list of bad edges, for each
element of the list there is a pointer from and to the corresponding 
edge in the main list of edges. This allows us to change the list of
bad edges after each step in constant time (although, we repeat, the
constant is impractically large).
\end{proof}

\subsection*{Outline of the computer search}
It remains to prove Lemma~\ref{noweight2}, and our proof of this
requires a computer. Unfortunately, both the number of labelings and
the number of possible cuts is far too large for a brute-force
approach: There are $2(2^9-1)-1$ edges of~$2T_9$, which means
more than $11^{1000}$ labelings, even if we use Lemma~\ref{noweight34}
to eliminate labeling with edges of weight 3 or~4. 
Moreover, there are roughly $(2^{2\cdot 2^8})^4$ internal cut labelings
in~$2T_9$, hence we cannot use brute-force even for one labeling.  To
overcome the second problem we shall recursively compute all of the
necessary information, a so-called ``menu'' for the subtrees, leading to
an efficient algorithm for a given labeling.  To solve the first
problem, instead of enumerating all labelings of~$2T_9$ and
computing the menu for them, we will iteratively find all menus
corresponding to all labelings of~$T_1$, $T_2$, \dots, $T_8$. This way
we avoid considering the same ``partial labeling'' several times.  To
further reduce the computational load, we will consider only ``worst
possible menus'' in each $T_i$.  Now, to the details.

\begin{figure}[ht]
\centerline{\includegraphics{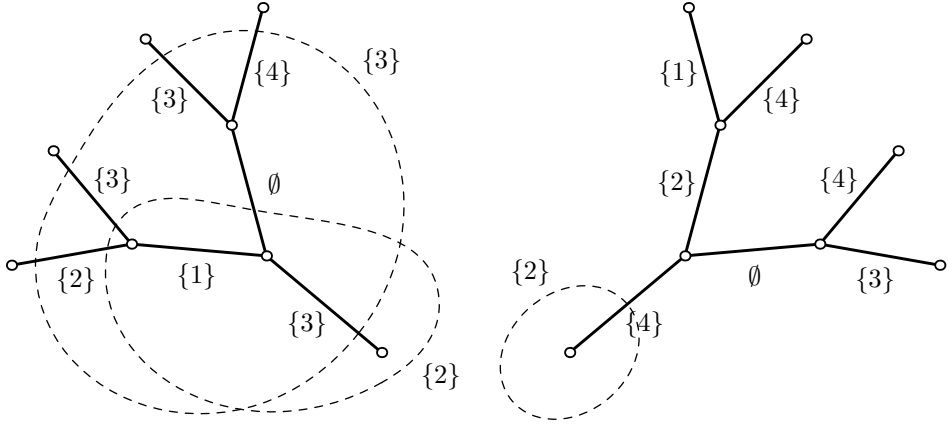}}
\caption{Illustration of the proof of Observation~\ref{parent_obs}:
We have $S=\emptyset$ and consider the 
internal $S$-swap indicated in Figure~\ref{fig:badlab}.
The split-up in this figure results in 
an internal $\{2\}$-swap in $T'$ (on the left) 
and an internal $\{2\}$-swap in $T''$ (on the right);
the set $Q$ equals $\{2\}$.}
\label{fig:splitup}
\end{figure}

If $S \subseteq [4]$ (we shall use $[4]$ to denote $\{1,2,3,4\}$), we define 
an internal cut labeling $Y$ of $T_i$ to be an \emph{internal $S$-swap} 
if $Y = (\delta(Z_1), \ldots, \delta(Z_4))$ where every 
$Z_i$ is a set of interior nodes (note that the root $r$ is an interior vertex) and 
$S = \{ i \in [4] : r \in Z_i \}$. Informally, 
an internal $S$-swap `switches~$S$ between the root and the leaves'
(see Figure~\ref{fig:splitup}).
A \emph{menu} is a mapping 
$M : \fsets \to \zet$.  If $T_i$ is a copy of a rooted tree with root $r$
and $X$ is a labeling of~$T_i$ then the \emph{menu corresponding
to~$X$} is defined as follows
\begin{equation}
\label{eq:menudef}
      M_X(S) = \min \{ \cost(X \sym Y) - \cost(X) : 
        \mbox{$Y$ is an internal $S$-swap} \} \,.
\end{equation}
Thus, the menu~$M_X$ associated with~$X$ is a function which tells us for
each subset $S \subseteq [4]$ the minimum cost of making an
internal $S$-swap. This is enough information to check whether we can
decrease the cost of a given labeling: if $T_1$, $T_2$, $T_3$
are trees meeting at a vertex and $X_i$~is the restriction of a 
labeling~$X$ to~$T_i$, then we can decrease the cost by a local swap 
(using only edges of $T_1$, $T_2$, and $T_3$)
if
we have $M_{X_1}(S) + M_{X_2}(S) + M_{X_3}(S) < 0$
for some~$S \in \fsets$.

For menus $M$, $N$ and a set $R \subseteq [4]$ we define 
$\Pm(M,N,R) : \fsets \rightarrow \zet$ to be the mapping given by the
following rule:
\begin{equation}\label{eq:menurec}
    \Pm(M,N,R)(S) = \min_{Q \in \fsets}
                  \Bigl( M(Q) + N(Q) 
                  +  a(|R \sym S \sym Q|) - a(|R|) \Bigr) \,.
\end{equation}
The motivation for this definition is the following observation, which
is the key to our recursive computation.  

\begin{observation} \label{parent_obs}
Let $X$ be a labeling of the tree $T_i$ ($i \ge 2$).  Let $e$ be the root edge of~$T_i$, 
let the two copies of $T_{i-1}$ that form $T_i - \{e\}$ be
denoted $T'$ and $T''$. Finally, let $X'$ and $X''$ be the restrictions
of the labeling $X$ to the trees $T'$ and $T''$.  Then 
$$
  M_X = \Pm(M_{X'},M_{X''},l_X(e)) \,.
$$
\end{observation}

\begin{proof}
Let $v$ be the end of the edge $e$ which is distinct from the root $r$.
Choose any $S \in \fsets$, we need to show, that
$M_X(S) = \Pm(M_{X'},M_{X''},l_X(e))(S)$,  
where the latter is defined by Equation~\eqref{eq:menurec}.

Consider an internal $S$-swap $Y = (\delta(Z_1), \ldots, \delta(Z_4))$ and
observe, that it is in 1-1 correspondence with a triple $(Y', Y'', Q)$, where
\begin{itemize}
  \item $Q = \{ i \in [4] : v \in Z_i \}$, 
  \item $Y' = (\delta_{T'}(Z_1 \cap V(T')), \ldots, \delta_{T'}(Z_4 \cap V(T')))$ 
     (here $\delta_{T'}$ means the neighborhood in $T'$), $Y'$ is an internal
     $Q$-swap in $T'$. Similarly 
  \item $Y'' = (\delta_{T''}(Z_1 \cap V(T'')), \ldots, \delta_{T''}(Z_4 \cap V(T'')))$
    is an internal $Q$-swap in $T''$. 
\end{itemize}
See also Figure~\ref{fig:splitup}, where the labeling of Figure~\ref{fig:badlab}
(described in detail below Observation~\ref{1cut}) is ``decomposed'' in this way.
With this correspondence we can decompose the change of cost between
labelings $X \sym Y$ and $X$ in the following way:

\begin{align*}
\cost (X \sym Y) - \cost (X) 
  =& \bigl(\cost (X' \sym Y') - \cost (X') \bigr) \\
   &+ \bigl(\cost (X'' \sym Y'') - \cost (X'') \bigr) \\
   &+ \bigl( a(| l_X(e) \sym S \sym Q |) - a(|l_X(e)|) \bigr) \\
\end{align*}

If we minimize the left-hand side over all internal $S$-swaps $Y$, we get~$M_X(S)$. 
Equivalently, we can minimize the right-hand side over all $Q \in \fsets$
and all internal $Q$-swaps $Y'$ (in $T'$) and $Y''$ (in $T''$).
However, for a fixed $Q$ the minimum over all such $Y'$ of
$\cost (X' \sym Y') - \cost (X')$ is $M_{X'}(Q)$, similarly for
the second summand. Thus, minimizing over $Q$, $Y$, and $Y'$ we get
the formula in Equation~\eqref{eq:menurec}.
\end{proof}

\def\calM{{\cal M}}
\def\calW{{\cal W}}

Using the above observation, it is relatively fast to compute the menu
associated with a fixed labeling of a tree $T_i$.  However, for our
problem, we need to consider all possible labelings of $T_i$.
Accordingly, we now define a few collections of menus which contain
all of the information we need to compute to resolve
Lemma~\ref{noweight2}.  Prior to defining these collections, we need
to introduce the following partial order on menus: if $M_1$ and $M_2$
are menus, we write $M_1 \lesseq M_2$ if $M_1(S) \le M_2(S)$ for every
$S \in \fsets$.  

We let $\calM_i$ be the set of all $M_X$, where $X$~is a labeling
of~$T_i$, and every $e \in E(T_i)$ satisfies $w_X(e) \le 2$.
We let $\calW_i$ denote the set of maximal (`worst') elements (with respect
to~$\lesseq$) of~$\calM_i$.  Further, we define two subsets of these sets: 
$\calM'_i$ denotes the set of menus corresponding to those
labelings~$X$ of~$T_i$ where each edge is of weight at most 2 and
where the root edge is labeled by~$\{1,2\}$. Finally, $\calW'_i$ is the
set of maximal elements of~$\calM'_i$.  The following observation
collects the important properties of these sets.

\begin{observation} \label{menus}
For every $i \ge 2$ we have
\begin{enumerate}[(1)]
\item
$
  \calM_i = \bigl \{ \Pm(M, N, R) 
          \mid M, N \in \calM_{i-1}, R \in \fsets, |R| \le 2 \bigr \}  
$
\item
$
  \calW_i = \max\limits_{\hbox{in $\lesseq$}}
      \bigl \{ \Pm(M, N, R) 
          \mid M, N \in \calW_{i-1}, R \in \fsets, |R| \le 2 \bigr \} 
$
\item
$
  \calW'_i = \max\limits_{\hbox{in $\lesseq$}}
      \bigl \{ \Pm(M, N, \{1,2\}) 
          \mid M, N \in \calW_{i-1} \bigr \}  
$
\end{enumerate}
\end{observation}

\begin{proof}
Part~(1) follows immediately from Observation~\ref{parent_obs}.  The
second part follows from this and from the fact that the mapping $\Pm$
is monotone with respect to the order $\lesseq$ on menus.  Part~(3)
follows by a similar argument.
\end{proof}

Next we state the key claim proved by our computer check.  

\begin{claim}[verified by computer]   \label{comp}
For every $W_1 \in \calW'_9$, and $W_2, W_3 \in \calW_8$ there exists 
$S \in \fsets$ such that $W_1(S) + W_2(S) + W_3(S) < 0$.
\end{claim}

We use Observation~\ref{menus} to give a practical scheme for
computing the collections $\calW_8$ and $\calW'_9$ followed by a simple
test for each possible triple.  Further details are described in the
Code Listing. With this, we are finally ready to give a proof of
Lemma~\ref{noweight2}.

\begin{proof}{\bf (of Lemma~\ref{noweight2})}
Let $X$ be an edge labeling of~$2T_9$ satisfying the assumptions; we
may suppose the central edge~$uv$ is labeled by~$\{1,2\}$.  Let $T^1$,
$T^2$, $T^3$ be the three distinct maximal subtrees of~$2T_9$ which
have $v$ as a leaf, and assume that $T^1$ contains the central
edge. Let $X_j$ denote the restriction of~$X$ to~$T^j$, and let
$M_j = M_{X_j}$ be the corresponding menu.
Choose $W_1 \in \calW'_9$, $W_2, W_3 \in \calW_8$ so
that $M_j \lesseq W_j$ holds for each $j$. By Claim~\ref{comp}, we may choose 
$S \in \fsets$ for which $W_1(S) + W_2(S) + W_3(S) <  0 $ and
by definition of~$\lesseq$ we have $M_1(S)+M_2(S)+M_3(S) < 0$, too.  
Let $X_j$ be the internal $S$-swap for which the minimum in the
definition of $M_j$ (Equation~(\ref{eq:menudef})) is attained.  Then
$Y = X_1 \sym X_2 \sym X_3$ is an internal cut labeling of $2T_9$ and 
$\cost(X \sym Y) - \cost(X) = M_1(S) + M_2(S) + M_3(S) < 0$.  
This completes the proof.
\end{proof}

\begin{remark}
In the definition of cost of a coloring, the values of parameters $a(i)$ can be
chosen in a variety of ways---provided we do penalize edges of weight~1.
Perhaps it seems more natural to have $a(1) = 0$, we only need to get
rid of the edges of weight $\ge 2$, so we might not penalize edges of weight~1
at all. However, this straightforward
approach does not work. Consider the edge labeling of~$2T_4$ the upper part of
which is depicted in Figure~\ref{fig:badlab}. (The lower part of $2T_4$ is a mirror 
image of this.) It is rather easy to verify, that
switching any local cut labeling does not get rid of edge of weight~2. Moreover,
this labeling can be extended to arbitrary~$2T_n$ by the `growing rules'
depicted in the lower part of the figure ($a$, $b$, $c$, $d$ stand for $\{1\}$, $\{2\}$, $\{3\}$,
$\{4\}$ in any order). 
On the other hand, by switching $\{2\}$ and $\{3\}$ on the cuts
depicted in the figure, we decrease the cost of the coloring by~$a(1)$.
Thus, choosing $a(1)$ nonzero allows us to distinguish, say, among 
various cut labelings where there is just one edge of weight~1. Then 
we can (by a series of local changes) move to a cut labeling, where we
can get rid of the edge of weight~1.

\begin{figure}[ht]
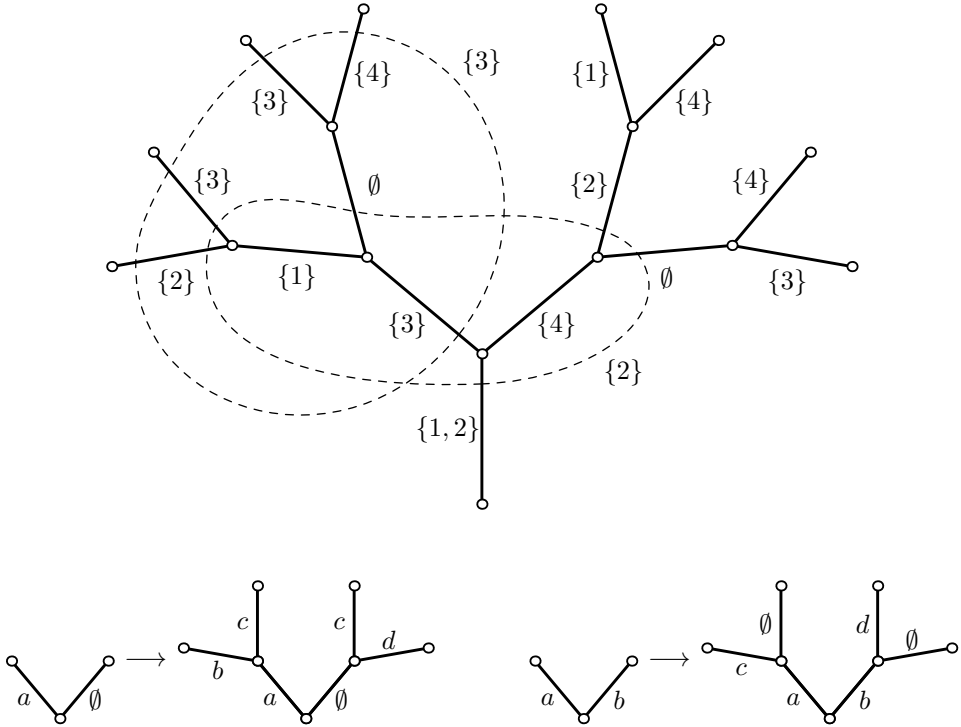

\centerline{\includegraphics{tree.0}}
\bigskip
\bigskip
\centerline{\includegraphics{tree.1} \raise 8mm \hbox{$\longrightarrow$}
  \includegraphics{tree.2}
  \hfill
\includegraphics{tree.3} \raise 8mm \hbox{$\longrightarrow$}
  \includegraphics{tree.4}}
\caption {A difficult labeling of $T_4$.}
\label{fig:badlab}
\end{figure}
\end{remark}

\begin{remark}
We note that we could prove Lemma~\ref{noweight34} by the same method
as Lemma~\ref{noweight2}; in fact a simple modification of the code
verifies both of these lemmas at the same time. The reason we put
Lemma~\ref{noweight34} separately is that it allows for an easy proof
by hand, and this hopefully makes the proof easier to understand. 

Another remark is that an easy modification of our method of
verifying Claim~\ref{comp} could decrease the running time by 30\%.
We did not want to obscure the main proof for this relatively small 
saving, but we wish to mention the trick here. In the process of 
enumerating the sets~$\calW_i$, we can throw away all menus~$M$
that satisfy $M(\emptyset) < 0$. It is not hard to show that
we still consider all `hard cases'.
\end{remark}

\begin{remark}
The necessity to use computer for huge amount of checking is
not entirely satisfying (although this point of view may be a rather
historically conditioned aesthetic criterion). It would be 
interesting to find a proof of Lemma~\ref{noweight2} without
extensive case-checking, perhaps by a careful inspection of the
sets~$\calW_i$.  
\end{remark}

\section{Some Equivalences}  \label{sec:equiv}

The goal of this section is to prove Proposition~\ref{equiv_prop} from
the Introduction (restated here for convenience as
Proposition~\ref{equiv_prop2}), which gives several graph properties
equivalent to the existence of a homomorphism to a projective cube
$PQ_{2k}$.  To prove this, it is convenient to first introduce another
family of graphs. For every positive integer~$n$, let $H_n$~denote
the graph with all binary vectors of length~$n$ forming the vertex set
and with two vertices being adjacent if they agree in exactly one
coordinate (note that $H_n$~is a Cayley graph on~$\zet_2^n$).

For odd $n$, the graph~$H_n$ has exactly two components, 
one containing all vertices with an even number of $1$'s, 
and the other all vertices with an odd number of $1$'s; 
we call the components $H_n^e$ and $H_n^o$, respectively.

\begin{observation} \label{pq_equiv}
For every $k \ge 1$ the graphs $H_{2k+1}^e$,  $H_{2k+1}^o$, and $PQ_{2k}$
are isomorphic.
\end{observation}

\begin{proof}
The mapping that sends each binary vector to its complementary vector gives an isomorphism 
between $H_{2k+1}^o$ and $H_{2k+1}^e$.  Thus, the simple graph obtained from $H_{2k+1}$ by 
identifying complementary vectors is isomorphic to $H_{2k+1}^e$ (and to
$H_{2k+1}^o$).  However, this graph is also isomorphic to $PQ_{2k}$,
since viewing the vertices of each as a pair of complementary vectors,
we see that $u$ and $v$ will be adjacent if and only if one vector
associated with $u$ and one vector associated with $v$ differ in
exactly 1 coordinate. 
\end{proof}

Now we are ready to prove the proposition.  

\begin{proposition} \label{equiv_prop2}
For every graph $G$ and nonnegative integer $k$, the following properties are equivalent.
\begin{enumerate}[(1)]
  \item There exist $2k$ pairwise disjoint cut complements.
  \item There exist $2k+1$ pairwise disjoint cut complements with union $E(G)$.
  \item $G$ has a homomorphism to $PQ_{2k}$.
  \item $G$ has a cut-continuous mapping to $C_{2k+1}$.  
\end{enumerate}
\end{proposition}

\begin{proof}
We shall show $(1) \implies (2) \implies (3) \implies (4) \implies (1)$.  

To see that $(1) \implies (2)$, let $S_1$, $S_2$, \ldots, $S_{2k}$ be
pairwise disjoint cut complements, and for every $1 \le i \le 2k$ let
$W_i = E(G) \setminus S_i$.  Now setting 
$S_{2k+1} = E(G) \setminus \cup_{1 \le i \le 2k} S_i =
E(G) \setminus \sym_{1 \le i \le 2k} W_i$ 
we have~(2).

Next we shall show that $(2) \implies (3)$.  
Let $S_1$, $S_2$, \ldots, $S_{2k+1}$ be $2k+1$ disjoint cut complements
with union $E(G)$ and for every $1 \le i \le 2k+1$ choose $U_i
\subseteq V(G)$ so that $S_i = E(G) \setminus \delta(U_i)$.  Now
assign to each vertex $v$ a binary vector $x^v$ of length $2k+1$ by
the rule $x^v_i = 1$ if $x \in U_i$ and $x^v_i = 0$ otherwise.  This
mapping gives a homomorphism from~$G$ to~$H_{2k+1}$, so by
Observation~\ref{pq_equiv} we conclude that $G$ has a homomorphism to~$PQ_{2k}$.  

Next we prove that $(3) \implies (4)$.  Since the composition of two
cut-continuous mappings is cut-continuous, it follows from
Observation~\ref{homocc} and Observation~\ref{pq_equiv} that it
suffices to find a cut-continuous mapping from $H_{2k+1}$ to $C_{2k+1}$.
To construct this, let $E(C_{2k+1}) = \{e_1,e_2,\ldots,e_{2k+1}\}$ and
define a mapping $g: E(H_{2k+1}) \rightarrow E(C_{2k+1})$ by the rule that
$g(uv) = e_i$ if $u$ and $v$ agree exactly in coordinate $i$.  We
claim that $g$ is a cut-continuous mapping.  To see this, let $R$ be a cut
of $C_{2k+1}$, let $J = \{ i \in \{1,2,\ldots,2k+1\} : e_i \in R \}$,
and note that $|J|$ is even.  Now let $X$ be the set of all binary
vectors with the property that there are an even number of $1$'s in
the coordinates specified by $J$.  Then $g^{-1}(R) = \delta(X)$ so our
mapping is cut-continuous as required.

To see that $(4) \implies (1)$, simply note that the preimage of any edge
of $C_{2k+1}$ is a cut complement, so the preimages of the $2k+1$
edges are $2k+1$ disjoint cut complements.  
\end{proof}

We can extract the key idea of the above proof as follows. 
Let $E_i \subseteq E(H_{2k+1})$ be the set of edges $uv$ such that
$u$ and~$v$ agree in exactly the $i$-th coordinate.\footnote{
If you think of~$H_n$ as of a Cayley graph, then $E_i$~consists of edges corresponding
to the $i$-th element of the generating set. We thank to Reza Naserasr
for this comment.}
The sets $E_1$, \dots, $E_{2k+1}$ form a partition of $E(H_{2k+1})$ into
disjoint cut complements.

\section{Code Listing} \label{sec:listing}

In this section we present the code used to verify Claim~\ref{comp}.
  The code is written in~C; it can be found at
\url{http://kam.mff.cuni.cz/~samal/papers/clebsch/} together with its 
output.
It runs about 30 minutes on a 2 GHz processor.\footnote{Over the course
  of the refereering process, this time decreased to 12 minutes on a 
  recent laptop.}
We have tested it with compilers gcc (version 3.0, 3.3, and 4.3), Intel~C, 
and Borland~C++ on several computers to minimize the possibility of error in the
proof due to erroneous computer hardware/software.

We use Observation~\ref{menus} to iteratively compute $\calW_{i+1}$
from $\calW_i$, this is accomplished by function~\code{W_update}.
By the same function we compute $\calW'_9$ from~$\calW_8$, 
we only provide a shorter (namely, one-element) list
of possible labels of the root edge.
Finally, we use~\code{final_test} to check whether all triples of menus
satisfy the inequality of Claim~\ref{comp}.
To simplify and speed up the code, we use static
data structures for $\calW_i$'s. That is, the elements of 
the set~$\calW_i$ are stored as \code{W[i][j]}, with 
$0 \le j < \code{W_size[i]}$ and with a limit \code{MAX=20000} on the
number of elements \code{W_size[i]}. 
If this number had turned out to be too small, the
program would have output an error message (this, however, did not happen).

Labels of edges, that is elements of $\fsets$ are represented as
integers from 0 up to 15.  For convenience variables that hold labels
have type \code{label} (which is a new name for \code{short}).
Symmetric difference of labels corresponds to bitwise xor---``\code{^}''.
Cost of edges are stored in variables of type \code{cost} 
(a new name for \code{int}). From Equation~\eqref{eq:menurec}
it is easy to deduce that $\Pm(M,N,R)(S) \le M(S) + N(S)$. 
Consequently, the largest coordinate of an element of~$\calW_i$
is in absolute value at most~$2^{i-1} a(4)$, and as we only use sets~$\calW_i$
for $i \le 9$, we will not have to store larger numbers
than an \code{int} can hold.
Other new data types are \code{menu} (array of 16 \code{cost}'s used to
represent a menu), and \code{comparison}---variables of that type are assigned
values $-1$, $0$, $1$, or \code{INCOMP=2} if the result of a corresponding
comparison (of two menus) is $\prec$, $=$, $\succ$ or incomparable.

When we need to compute $M = \Pm(M_1, M_2, c)$, this is implemented as 
\code{add_menus(M_1,M_2,children); p_menu(children, parent, M)}.
(The reason for this two-step process is that \code{children} is only computed
once and then used for all possible $c$'s.)
Here \code{children} corresponds to the sum $M_1 + M_2$, \code{parent} is a menu
corresponding to the single edge of~$T_1$ labeled by $c$. 
Then we insert the menu in the set~$\calW_i$ (array \code{W[i]})
by calling \code{insert_menu}. 
This simply compares $M$ to all menus in \code{W[i]}. If some of them is
$\succ M$, we are done with $M$. Otherwise, we add $M$ to \code{W[i]} and
delete all menus in \code{W[i]}, that are possibly $\prec M$. 
(This is implemented in a somewhat roundabout way (to save time).
To fill the empty spaces after the deleted menus we move there menus
from the end, that is \code{W[i][W_size[i]-1]}. This avoids moving
all of the menus in memory.
When we implemented the deletion of `small' menus in this function 
in a more straightforward manner (`move everything left'), the running 
time did approximately double.)

{ 
\begin{lstlisting}
#include <stdio.h>
#include <limits.h>
#define MAX 20000    // limit on size of the sets W_i

typedef short label;  
typedef int cost;
typedef cost menu[16];

typedef short comparison;
comparison INCOMP = 2;  

cost a[5]={0,1,10,40,1000};
cost labelcost[16];   // cost of edge labeled by each possible label
menu one_label[1];    // W'_1, i.e. one_label[0] corresponds to T1 labeled by {1,2}
menu W[9][MAX];
menu Wprime[MAX];     // W'_9
int Wsize[9];         // Wsize[i] is the number of elements of  W[i]
int Wprimesize;       // the number of elements of Wprime

void menu_from_label(label r, menu M) {
// M will be the menu corresponding to T1 labeled by r
  label s;
  for(s=0; s<16; s++) 
      M[s] = labelcost[r ^ s] - labelcost[r];
}

void init_variables() {
  label s;
  for (s=0; s<16; s++)
    labelcost[s] = a[(s&1) + ((s>>1) & 1) + ((s>>2) & 1) + ((s>>3)&1)]; 
// the right hand side is a[n], where n is the number of ones 
// in binary representation of s

  menu_from_label(3,one_label[0]); // 3 corresponds to {1,2}

  Wsize[1]=0;
  for(s=0; s<16; s++) 
    if (labelcost[s] < a[3]) 
      menu_from_label(s,W[1][Wsize[1]++]);
}

void add_menus(menu M1, menu M2, menu sum) {
  label s;
  for (s=0; s<16; s++) 
    sum[s] = M1[s]+M2[s];
}

comparison sign(int n) {
  if (n > 0) return 1;
  if (n < 0) return -1;
  return 0;
}

comparison compare_menus(menu M1, menu M2) {
// returns  -1, 0, 1, INCOMP, depending on
// whether M1<M2, M1=M2, M1 > M2, or they are incomparable
  label s;
  comparison t, current=0;

  for (s=0; s < 16; s++) {
    t = sign (M1[s] - M2[s]);
    if ((t != 0) && (t == -current)) return INCOMP;
    if (current == 0) current = t;
  }
  return current;
}

void p_menu(menu children, menu parent, menu output) {
// children is the sum of the menus of the two subtrees
// parent corresponds to the root edge
  label s, q;
  cost new, current_best;

  for (s=0; s<16; s++) {
    current_best = children[0] + parent[s];  // for q=0
    for (q = 1; q < 16; q++) {
      new = children[q] + parent[s ^ q];     // using equation (2)
      if (new < current_best) current_best = new;
    }
    output[s] = current_best;
  }
}

void insert_menu(menu *book, int *booksize, menu M) {
// book is an array of menus: book[0] ... book[*booksize-1]
// booksize is the number of elements of book, we are inserting M
  int i;
  label s;
  comparison t=0;

  for (i=0; i < *booksize; i++) {
    t = compare_menus(M,book[i]);
    if (t <= 0) return;  // M <= book[i], so we will not insert M
    if (t == 1) break;   // M > book[i], so no other element
      // of book may be larger than M
  }

// either M is INCOMP with every menu
// or: 
  if (t==1)  // i.e. M > book[i], we will 
      // delete all elements of book that are < M
    for ( ; i < *booksize; i++) {
      while (i <  *booksize && compare_menus(M,book[i])==INCOMP) 
        i++;
      // at this point we have found another menu less than M, namely book[i]
      while (*booksize > i && compare_menus(book[*booksize-1],M) <= 0 ) 
        (*booksize)--; // we abandon small menus at the end
      if (*booksize <= i)  // all the remaining menus were small
        break;         // all menus < M are deleted
      // there is a big menu at the end, we move it to book[i]:
      (*booksize)--;
      for (s = 0; s<16; s++) 
        book[i][s] = book[*booksize][s];
    }

// we insert M as the last element of book
  if (*booksize == MAX) printf("too short array!\n");
  else {
    for (s = 0; s<16; s++)
      book[*booksize][s] = M[s];
    (*booksize)++;
  }
}

void W_update(menu *oldW, int oldsize, menu *root_edge, int rootsize, 
              menu *newW, int *newsize) {
  menu N, children;
  int i, j, k;

  *newsize = 0;
  for (i=0; i < oldsize; i++)
    for (j=i; j < oldsize; j++) {
      add_menus(oldW[i],oldW[j],children);
      for (k=0; k < rootsize; k++) {
        p_menu(children,root_edge[k],N);
        insert_menu(newW, newsize, N);
      }
    }
}


int final_test(menu *C, int Csize, menu *P, int Psize) {
  int i, j, k;
  label s;
  int counter=0; // number of found counterexamples to Claim 2.7
  menu children;

  for (i=0; i < Csize; i++)
    for (j=i; j < Csize; j++) {
      add_menus(C[i],C[j],children);
      for (k=0; k < Psize; k++) {
        counter ++; // we are testing possible counterexample P[k], C[i], C[j]
        for (s=0; s<16; s++)
          if (children[s]+P[k][s] < 0) { counter--; break;}
        // Claim 2.7 holds for P[k],C[i],C[j]
        // we proceed by testing another triple
      }
    }
  return counter;
}

int main() {
  int i;
  init_variables();

  printf("%d %d\n", INT_MIN, INT_MAX); 
   // to check whether 2^8.1000 is not too large
  for (i=1; i<8; i++) {
    W_update(W[i],Wsize[i],W[1],Wsize[1],W[i+1],&Wsize[i+1]);
    printf("The size of W_%d is: %d\n",i+1,Wsize[i+1]);
  }
  W_update(W[8],Wsize[8],one_label,1,Wprime,&Wprimesize);

  if (final_test(W[8],Wsize[8],Wprime,Wprimesize) == 0) 
    printf("\nProof is finished.\n\n"); 

  return 0;
}
\end{lstlisting} }

\section* {Acknowledgements}

We thank to Stephen C. Locke, Nicholas C. Wormald and Brendan McKay for helpful
discussion. We thank to Daniel Kr\'al' for pointing to us the coloring in
Figure~\ref{fig:badlab}.

\bibliographystyle{rs-amsplain}
\bibliography{clebsch}

\providecommand{\bysame}{\leavevmode\hbox to3em{\hrulefill}\thinspace}
\providecommand{\MR}{\relax\ifhmode\unskip\space\fi MR }
\providecommand{\MRhref}[2]{%
  \href{http://www.ams.org/mathscinet-getitem?mr=#1}{#2}
}
\providecommand{\href}[2]{#2}
\begin{thebibliography}{10}

\bibitem{Biggs-cubic}
Norman Biggs, \emph{Constructions for cubic graphs with large girth}, Electron.
  J. Combin. \textbf{5} (1998), Article 1, 25 pp. (electronic).

\bibitem{BL}
J.~Adrian Bondy and Stephen~C. Locke, \emph{Largest bipartite subgraphs in
  triangle-free graphs with maximum degree three}, J. Graph Theory \textbf{10}
  (1986), no.~4, 477--504.

\bibitem{DNR}
Matt DeVos, Jaroslav Ne{\v s}et{\v r}il, and Andr{\'e} Raspaud, \emph{On
  edge-maps whose inverse preverses flows and tensions}, Graph Theory in Paris:
  Proceedings of a Conference in Memory of Claude Berge (J.~A. Bondy,
  J.~Fonlupt, J.-L. Fouquet, J.-C. Fournier, and J.~L.~Ramirez Alfonsin, eds.),
  Trends in Mathematics, Birkh{\"a}user, 2006, pp.~109--138.

\bibitem{GL-unique}
Don Greenwell and L{\'a}szl{\'o} Lov{\'a}sz, \emph{Applications of product
  colouring}, Acta Math. Acad. Sci. Hungar. \textbf{25} (1974), 335--340.

\bibitem{Guenin-joins}
Bertrand Guenin, \emph{Packing {T}-joins and edge colouring in planar graphs},
  (to appear).

\bibitem{HH}
Roland H{\"a}ggkvist and Pavol Hell, \emph{Universality of {$A$}-mote graphs},
  European J. Combin. \textbf{14} (1993), no.~1, 23--27.

\bibitem{Hatami-C7}
Hamed Hatami, \emph{Random cubic graphs are not homomorphic to the cycle of
  size~7}, J. Combin. Theory Ser. B \textbf{93} (2005), no.~2, 319--325.

\bibitem{HatamiZhu}
Hatami Hatami and Xuding Zhu, \emph{The fractional chromatic number of graphs
  of maximum degree at most three}, submitted.

\bibitem{HT-independenceratio}
Christopher~Carl Heckman and Robin Thomas, \emph{A new proof of the
  independence ratio of triangle-free cubic graphs}, Discrete Math.
  \textbf{233} (2001), no.~1-3, 233--237, Graph theory (Prague, 1998).

\bibitem{Hladky}
Jan Hladk\'{y}, \emph{Bipartite subgraphs in a random cubic graph}, Bc.~thesis,
  Charles University, 2006, http://kam.mff.cuni.cz/\~{}hladky/bak.pdf.

\bibitem{HS}
Glenn Hopkins and William Staton, \emph{Extremal bipartite subgraphs of cubic
  triangle-free graphs}, J. Graph Theory \textbf{6} (1982), no.~2, 115--121.

\bibitem{KNS}
Alexandr~V. Kostochka, Jaroslav Ne{\v s}et{\v r}il, and Petra
  Smol{\'\i{}}kov{\'a}, \emph{Colorings and homomorphisms of degenerate and
  bounded degree graphs}, Discrete Math. \textbf{233} (2001), no.~1-3,
  257--276, Fifth Czech-Slovak International Symposium on Combinatorics, Graph
  Theory, Algorithms and Applications, (Prague, 1998).

\bibitem{McK-proc}
Brendan McKay, \emph{Maximum bipartite subgraphs of regular graphs with large
  grith}, Proceedings of the 13th Southeastern Conf. on Combinatorics, Graph
  Theory and Computing, Boca Raton, Florida, 1982.

\bibitem{Naserasr-thesis}
Reza Naserasr, \emph{Homomorphisms and edge colourings of planar graphs}, Ph.D.
  thesis, Simon Fraser University, 2003.

\bibitem{Nes-Aspects}
Jaroslav Ne{\v{s}}et{\v{r}}il, \emph{Aspects of structural combinatorics (graph
  homomorphisms and their use)}, Taiwanese J. Math. \textbf{3} (1999), no.~4,
  381--423.

\bibitem{dMN}
Jaroslav Ne{\v{s}}et{\v{r}}il and Patrice Ossona~de Mendez, \emph{Grad and
  classes with bounded expansion. {III}. {R}estricted graph homomorphism
  dualities}, European J. Combin. \textbf{29} (2008), no.~4, 1012--1024.

\bibitem{rs-thesis}
Robert {\v S}\'amal, \emph{On {XY} mappings}, Ph.D. thesis, Charles University,
  2006.

\bibitem{NS-TT1}
Jaroslav~Ne\v set\v ril and Robert {\v S}\'amal, \emph{Tension-continuous
  maps---their structure and applications}, submitted, arXiv:math.CO/0503360.

\bibitem{NS-TT2}
Jaroslav~Ne\v set\v ril and Robert {\v S}\'amal, \emph{On tension-continuous
  mappings}, European J. Combin. \textbf{29} (2008), no.~4, 1025--1054,
  arXiv:math.CO/0602563.

\bibitem{seymour-r}
Paul~D. Seymour, \emph{On multicolourings of cubic graphs, and conjectures of
  {F}ulkerson and {T}utte}, Proc. London Math. Soc. (3) \textbf{38} (1979),
  no.~3, 423--460.

\bibitem{WW}
Ian~M. Wanless and Nicholas~C. Wormald, \emph{Regular graphs with no
  homomorphisms onto cycles}, J. Combin. Theory Ser. B \textbf{82} (2001),
  no.~1, 155--160.

\bibitem{Wormald-survey}
Nicholas~C. Wormald, \emph{Models of random regular graphs}, Surveys in
  combinatorics, London Math. Soc. Lecture Note Ser., vol. 267, Cambridge Univ.
  Press, Cambridge, 1999, pp.~239--298.

\bibitem{Zy}
Ond{\v r}ej Z{\'y}ka, \emph{On the bipartite density of regular graphs with
  large girth}, J. Graph Theory \textbf{14} (1990), no.~6, 631--634.

\end{thebibliography}

\end{document}